\documentclass[11pt]{amsart}

\usepackage{graphics,graphicx}
\usepackage{tikz-cd}
\usetikzlibrary{calc}

\usepackage{dynkin-diagrams}

\usepackage{dsfont}

\usepackage{xcolor}

\usepackage{enumitem}

\usepackage{booktabs}

\usepackage[a4paper]{geometry}
\usepackage[utf8]{inputenc}
\usepackage[T1]{fontenc}

\usepackage{url}

\usepackage{amssymb,amsmath,amsthm}

\usepackage[backref=page,colorlinks=true,allcolors=blue]{hyperref} %

\newtheorem{thm}{Theorem}[section]

\newtheorem{conj}[thm]{Conjecture}
\newtheorem{pbm}[thm]{Problem}

\newtheorem*{ack}{Acknowledgements}

\theoremstyle{definition}
\newtheorem{defn}[thm]{Definition}

\newtheorem{exa}[thm]{Example}

\DeclareMathOperator{\SL}{SL}
\DeclareMathOperator{\SO}{SO}

\DeclareMathOperator{\Stab}{Stab}

\DeclareMathOperator{\relint}{RelInt}

\DeclareMathOperator{\Ric}{Ric}
\DeclareMathOperator{\Bl}{Bl}

\DeclareMathOperator{\Hom}{Hom}

\DeclareMathOperator{\DF}{DF}
\DeclareMathOperator{\MNA}{M^{NA}}
\DeclareMathOperator{\nld}{nld}

\newcommand{\bbC}{\mathbb{C}}
\newcommand{\bbQ}{\mathbb{Q}}
\newcommand{\bbZ}{\mathbb{Z}}
\newcommand{\bbP}{\mathbb{P}}
\newcommand{\bbR}{\mathbb{R}}
\newcommand{\bbN}{\mathbb{N}}

\begin{document}

\title{On the effective YTD conjecture}

\author{Thibaut Delcroix}
\address{Thibaut Delcroix, Univ Montpellier, CNRS, Montpellier, France}
\email{thibaut.delcroix@umontpellier.fr}
\urladdr{http://delcroix.perso.math.cnrs.fr/}

\date{\today}

\begin{abstract}
We formulate an effective variant of the Yau-Tian-Donaldson conjecture, then review effective results on K-stability of spherical varieties, that is, K-stability criterions which can be effectively computed given the combinatorial data associated with the variety. We focus on the standard notion of K-stability as defined by Donaldson for constant scalar curvature Kähler metrics. 
\end{abstract}

\maketitle


The quest for canonical Kähler metrics has been a leading direction of research in complex geometry since Calabi's program proposed in his lecture at the 1954 International Congress of Mathematicians. 
According to the general philosophy in complex geometry that differential and algebraic aspects are strongly intertwined, and inspired both by GIT pictures and Kobayashi-Hitchin correspondence, a conjecture relating existence of canonical Kähler metrics with an algebro-geometric notion of K-stability quickly appeared. It is now known as the Yau-Tian-Donaldson, or YTD conjecture.

The last decade has seen tremendous progress around this conjecture, especially in the Fano case which is completely settled. 
Even more interestingly, in the Fano case, the resolution of the YTD conjecture has allowed a huge jump in the possibilities of verifying effectively the existence of Kähler-Einstein metrics on a given Fano manifold. 
In this short survey, we explore the possibility of similar applications of the resolution of the YTD conjecture in the case of constant scalar curvature Kähler metrics, focusing on the case of spherical manifolds, a large class of manifolds with a large group of automorphisms. 

\begin{ack}
The author is partially funded by ANR-21-CE40-0011 JCJC project MARGE.
\end{ack}

\section{K-stability and the Yau-Tian-Donaldson conjectures}

\subsection{CscK metrics}

Let \(X\) be a complex manifold, and let \(\omega\) be a Kähler form on \(X\), referred to as a Kähler metric as usual in the field. 
In particular, it defines a Riemannian metric on \(X\), with associated notions of curvature. 
In the Kähler world, the Ricci curvature may be very conveniently encoded by a closed \((1,1)\)-form called the \emph{Ricci form}. 
In local holomorphic coordinates \((z_1,\ldots,z_n)\) for the complex manifold \(X\), the Kähler form \(\omega\) is defined by a local potential \(\varphi\): 
\[ \omega =i\partial\bar\partial \varphi = \sum_{j,k=1}^n \frac{\partial^2\varphi}{\partial z_j\partial \bar{z}_k} idz_j\wedge d\bar{z}_k \]
and the Ricci form is then given by 
\[ \Ric(\omega) = i\partial\bar\partial \left(-\ln\det \left(\frac{\partial^2\varphi}{\partial z_j\partial \bar{z}_k}\right) \right) \]
A \emph{Kähler-Einstein metric} is a Kähler metric \(\omega\) with constant Ricci curvature, that is, a solution to 
\[ \Ric(\omega)= t\omega\]
for some real number \(t\). 
This is the best known notion of canonical Kähler metrics, but unfortunately may only exist on very particular manifolds since the above equation implies a strong condition on the first Chern class: \(c_1(X)\) is either Kähler, anti-Kähler or zero. 
To solve this issue, one considers a weaker notion of constant curvature. 

\begin{defn}
A Kähler form \(\omega\) on a complex manifold \(X\) is \emph{cscK} if its scalar curvature \(S_{\omega}:=n\frac{\Ric(\omega)\wedge \omega^{n-1}}{\omega^n}\) is a constant function. 
\end{defn}

If \(\omega\) is cscK, then the value of its (constant) scalar curvature is encoded by its Kähler class. 
Furthermore, if there exists a cscK metric in a Kähler class, then it is unique up to the action of the automorphism group \cite{Berman-Berndtsson_2017,Chen-Tian_2008,Donaldson_2001}. 
As a consequence, it is a canonical metric for the data of a compact complex manifold \(X\) \emph{and} a fixed Kähler class \(\alpha\). 
For this survey, we will focus on the case when the Kähler class is that of an ample line bundle for simplicity, although some of the discussion carries over to the transcendental case, see in particular \cite{Dervan-Ross_2017,Sjostrom_Dyrefelt_2020} for the definition of K-stability for a transcendental Kähler class. 

\subsection{K-stability}

The definition used for K-stability in this paper was introduced by Donaldson in 2002 \cite{Donaldson_2002}. 
Several variants were introduced throughout the years, we will use the relations between some of these as discussed in \cite{Boucksom-Hisamoto-Jonsson_2017}. 

\begin{defn}
A \emph{test configuration} for \((X,L)\) is a normal variety \(\mathcal{X}\) equipped with a projective flat morphism \(\pi : \mathcal{X}\to \bbC\), a \(\pi\)-ample line bundle \(\mathcal{L}\) and an \(\mathcal{L}\)-linearized \(\bbC^*\)-action such that \(\pi\) is \(\bbC^*\)-equivariant with respect to the action by multiplication on \(\bbC\), and the fiber \((\pi^{-1}(1),\mathcal{L}|_{\pi^{-1}(1)})\) is isomorphic to \((X, L^r)\) for some positive integer \(r\).  
\end{defn}

In the following, \(\mathcal{X}\) denotes a test configuration for \((X,L)\) (the line bundle, the action and the morphism are omitted from the notation). 
The most important data in a test configuration is that of its central fiber \(\mathcal{X}_0=\pi^{-1}(0)\) as a scheme-theoretic fiber, which may in general be non-reduced and have several components. 
It is equipped with an action of \(\bbC^*\) since \(\pi\) is equivariant. As a consequence one has a family of representations of \(\bbC^*\) given by the spaces 
\[ H^0(\mathcal{X}_0,\mathcal{L}|_{\mathcal{X}_0}^k) \]
for \(k\in \bbZ_{>0}\). 
Let \(d_k\) denote the dimension of this space, and let \(w_k\) denote the total weight of this representation. 
Then one has an expansion 
\[ \frac{w_k}{kd_k} = F_0-F_1\frac{1}{k}+o\left(\frac{1}{k}\right) \]

\begin{defn}
The \emph{Donaldson-Futaki invariant} \(\DF(\mathcal{X})\) of the test configuration \(\mathcal{X}\) is the coefficient \(F_1\) in this expansion. 
\end{defn}

It turns out that the Donaldson-Futaki invariant defined this way does not behave well under the base changes of the form \(\bbC\to \bbC, z\mapsto z^d\), which can always be used to transform a test configuration to a test configuration with reduced central fiber. 
A variant called the non-Archimedean Mabuchi functional is enough to check K-stability, it coincides with the Donaldson-Futaki invariant on test configurations with reduced central fibers, and is linear in \(d\) along base changes as above.
We denote by \(\MNA(\mathcal{X})\) the non-Archimedean functional applied to the test configuration \(\mathcal{X}\). 

There are some remarkable classes of test configurations, depending on properties of their central fiber. 

\begin{defn}
A test configuration \(\mathcal{X}\) is a \emph{special test configuration} if the scheme-theoretic fiber \(\mathcal{X}_0\) is a normal variety. 
It is a \emph{product test configuration} if \(\mathcal{X}_0\) is isomorphic to \(X\). 
If \(X\) is equipped with a holomorphic action of a reductive group \(G\), then a test configuration \(\mathcal{X}\) is \emph{\(G\)-equivariant} if the polarized total space \((\mathcal{X},\mathcal{L})\) is equipped with a \(G\)-action which commutes with the \(\bbC^*\)-action, and the isomorphism between \(\pi^{1}\) and \(X\) is \(G\)-equivariant. 
\end{defn}

We may finally define the notion of K-polystability introduced by Donaldson in \cite{Donaldson_2002}.

\begin{defn}
A polarized variety \((X,L)\) is \emph{K-polystable} if, for any test configuration \(\mathcal{X}\) for \((X,L)\), \(\DF(\mathcal{X})\geq 0\), and equality holds only if  \(\mathcal{X}\) is a product test configuration. 
\end{defn}

In view of the above remarks, one may replace \(\DF\) by \(\MNA\) in the definition, and one may consider only test configurations with reduced central fiber.

\subsection{Effective Yau-Tian-Donaldson conjecture for cscK metrics}

A first precise statement of the Yau-Tian-Donaldson conjecture for cscK metrics follows from Donaldson's definition of K-stability. 

\begin{conj}
\label{YTD}
Let \((X,L)\) be a smooth polarized manifold. 
There exists a cscK metric in \(c_1(L)\) if and only if \((X,L)\) is K-polystable. 
\end{conj}

This conjecture has been proved for \(L=K_X^{-1}\) (in particular, \(X\) is a Fano manifold) by Chen-Donaldson-Sun \cite{CDS1,CDS2,CDS3} and Tian \cite{Tian_2015}. 
In fact, a cscK metric in \(c_1(X)=c_1(K_X^{-1})\) is a Kähler-Einstein metric, so the result provides a completely algebro-geometric characterization of existence of Kähler-Einstein metrics, together with Aubin and Yau's foundational existence results \cite{Aubin_1978,Yau_1978} for \(c_1(X)\leq 0\). 
Actually, a much better result is proved, showing that it is enough to consider only special test configurations (this was much closer to the original version of K-stability introduced by Tian \cite{Tian_1997}, and from the K-stability point of view, this was shown earlier for \((X,K_X^{-1})\) a \(\bbQ\)-Fano variety by Li and Xu \cite{Li-Xu_2014}).  
Further refinements of this result, and new proofs were obtained, allowing in particular to reduce to test configurations that are equivariant with respect to a reductive group of automorphism. 

\begin{thm}[\cite{Datar-Szekelyhidi_2016}]
Let \(X\) be a Fano manifold equipped with an action of a reductive group \(G\). Then there exists a Kähler-Einstein metric on \(X\) if and only if \(\DF(\mathcal{X})\geq 0\) for any special \(G\)-equivariant test configuration for \((X,K_X^{-1})\), with equality if and only if the test configuration is product.  
\end{thm}

A further motivating element for this type of result was the Hamilton-Tian conjecture originating from \cite{Tian_1997} as well, claiming that the Gromov-Hausdorf limit of the Kähler-Ricci flow on a Fano manifold should be a Kähler-Ricci soliton, possibly on a degeneration of \(X\). 
This conjecture was proved in several papers \cite{Chen-Wang_2012,Bamler_2018,Wang-Zhu_2021}, and the degeneration of \(X\) on which the Kähler-Ricci soliton lives is the central fiber of a (two-step) special test configuration for \((X,K_X^{-1})\) by \cite{Chen-Sun-Wang_2018}. 

We suggest a (very optimistic) generalization of this reduction to special test configurations in the Fano case, called the effective YTD conjecture. 
Here there is no known analogue of the Kähler-Ricci flow to suggest the nature of the destabilizing objects (although the Calabi flow is a candidate). 

\begin{pbm}[Effective YTD conjecture]
\label{effective_YTD}
Let \(n\) and \(\rho\) be positive integers, prove that there exists a positive integer \(m=m(n,\rho)\) such that: for any smooth projective manifold \(X\) of dimension \(n\) and Picard rank \(\rho\), and any ample line bundle \(L\) on \(X\),  
there exists a cscK metric in \(c_1(L)\) if and only if \(\DF(\mathcal{X})\geq 0\) for any test configuration, equivariant under a Levi subgroup of its reduced automorphism group, whose central fiber is reduced with at most \(m\) irreducible components, and equality holds if and only if the test configuration is product.  
Estimate the value of \(m\). 
\end{pbm}

We will highlight, in the following, positive results in this direction, that will also explain why we call this an \emph{effective} YTD conjecture. 
Verifying the conjecture for other classes of polarized manifolds seems to be the best hope towards the conjecture currently: it is for example completely open for surfaces, as well as for three-dimensional toric polarized manifolds.  

Although there are, to the author's knowledge, no known counterexample to Conjecture~\ref{YTD}, it was soon realized that a stronger notion of K-stability was likely needed to prove a version of the YTD conjecture for cscK metrics (in particular, with the examples of \cite{ACGTIII_2008}).
For this one can either enlarge the space of conditions from test configurations to real test configurations or filtrations or to objects in the finite energy completion of test configurations, or one can consider a uniform version of K-stability. 
Recent results provide a full solution of the YTD conjecture for this strengthened notion of K-stability \cite{Li_2022,BJ}.
While these results are the culmination of a series of works towards the YTD conjecture, the notion of K-stability involved seems in turn to be further away from an effectively checkable condition. 
It leaves fully open the algebro-geometric side of Problem~\ref{effective_YTD}, which precisely aims at reducing the space of conditions to check K-stability.  
On the other hand, one may readily reformulate Problem~\ref{effective_YTD} with real test configurations whose central fibers have a controlled number of irreducible components without losing the effective aspect. 

\section{Spherical varieties}

As highlighted above, it is desirable, and proved in the various versions of the YTD conjecture known to hold, to be able to reduce to test configurations invariant under the action of a reductive group \(G\) acting holomorphically on \(X\). 
If the group \(G\) is a torus, then the largest possible effective action is that of a torus with the same dimension as that of \(X\), in which case the manifold \(X\) is a toric manifold. 
These toric manifolds were instrumental in the development of the theory \cite{Wang-Zhu_2004,Donaldson_2002}. 
We now focus on a wide generalization of toric manifolds, the class of spherical manifolds. 

\subsection{The multiplicity free property}

Let \(G\) be a connected complex linear reductive group. 
We fix a maximal torus \(T\) of \(G\), a Borel subgroup \(B\) of \(G\) containing \(T\), and we let \(\Phi\) and \(\Phi^+\) denote the corresponding roots and positive roots of \(G\). 
Any finite dimensional representation \(V\) of \(G\) is a direct sum of irreducible representations 
\(V\simeq \bigoplus V_{\lambda}^{m_{\lambda}}\) 
where \(\lambda\) runs over dominant weights of \(G\), \(V_{\lambda}\) is an irreducible representation of \(G\) with highest weight \(\lambda\), and \(m_{\lambda}\) is called the multiplicity. 
The representation \(V\) is called \emph{multiplicity free} if for all \(\lambda\), \(m_{\lambda}\in \{0,1\}\). 

\begin{defn}
A projective normal \(G\)-variety \(X\) is \emph{spherical} or \emph{multiplicity free} if for any \(G\)-linearized line bundle \(L\) on \(X\), the \(G\)-representation \(H^0(X,L)\) is multiplicity free. 
\end{defn}

The above definition is neither the most usual definition nor the easiest to check on examples, but highlights why spherical varieties are somewhat easier to deal with in view of the definition of K-stability. 
The most usual definition is that \(B\) acts with a Zariski open orbit on \(X\), and the proof of the equivalence of the two definition may be found in \cite[Theorem~2.1.2]{Perrin_2014}.

\begin{exa}
A generalized flag variety \(G/P\) (that is, the isotropy group \(P\) is a parabolic subgroup of \(G\) or, equivalently, the homogeneous space \(G/P\) is projective) is spherical under the action of \(G\). 
In this case, for any \(G\)-linearized line bundle \(L\) on \(G/P\), \(H^0(X,L)\) is an irreducible representation of \(G\). 
\end{exa}

\begin{exa}
A toric variety \(X\), that is, a normal variety equipped with an effective action of a \(\dim(X)\)-dimensional torus \(T\), is spherical under the action of \(T\). 
In this case the space of sections of an ample \(T\)-linearized line bundle \(L\) on \(G/P\) decomposes as a sum of irreducible representations, thus of dimension one, each encoded by the character through which \(T\) acts. 
Furthermore, the characters actually appearing are given by the integral points of a polytope in the vector space generated by characters. 
We will see how this description generalizes for spherical varieties in the next section.  
\end{exa}

One can put together the two classes to obtain further examples of spherical varieties. 

\begin{defn}
Assume that \(X\) is a normal \(G\)-variety, and \(x\in X\) is such that \(G\cdot x\) is Zariski open in \(X\) and \(\Stab_G(x)\) contains the unipotent radical of a Borel subgroup of \(G\). 
Then \(X\) is a spherical \(G\)-variety, called \emph{horospherical}. 
\end{defn}

In that case the open orbit \(G\cdot x\) is a torus bundle over a generalized flag variety \(G/P\), and the closure in \(X\) of the fiber over any point \(b\in G/P\) is a toric variety under the action of (a quotient of) the center of a Levi subgroup of \(P\). 

\begin{exa}
For any \(k\in \bbZ_{\geq 1}\), the Hirzebruch surface \(\bbP_{\bbP^1}(\mathcal{O}\oplus \mathcal{O}(k))\), equipped with the natural lift of the action of \(\SL_2\) on \(\bbP^1\), is a horospherical varieties. 
The weighted projective space \(\bbP(1,1,k)\) obtained by contracting the negative section is horospherical as well. 
\end{exa}

This class contains varieties that are very different from both generalized flag varieties and toric varieties. 
For example, Boris Pasquier \cite{Pasquier_2009} classified among these horospherical varieties the smooth, Picard rank one examples and obtained several infinite families which are not homogeneous. 
By contrast, the only smooth, Picard rank one toric manifolds are the projective spaces. 

Furthermore, spherical varieties contain many subclasses which are far from horospherical varieties, perhaps the most famous of which being De Concini and Procesi's wonderful compactifications of complex symmetric spaces \cite{De_Concini-Procesi_1983}. 
Recently, together with Pierre-Louis Montagard we classified spherical actions of low rank (see the next section for the definition of the rank) on Fano manifolds of low dimensions \cite{Delcroix-Montagard}. 
For example the nine Fano threefolds with spherical actions which admit neither a generalized flag manifold nor a toric structure are the Fano threefolds number 2-29, 2-30, 2-31, 3-19, 3-20, 3-22, 3-24, 4-7, 4-8 in Mori-Mukai numbering \cite{Mori-Mukai_1981}. 

\begin{exa}
The blowup \(\Bl_{\mathcal{Q}^1}(\mathcal{Q}^3)\) of the three-dimensional quadric \(\mathcal{Q}^3\) along a conic \(\mathcal{Q}^1\) is a spherical variety under the action of \(\SL_2(\bbC)\times \bbC^*\simeq \SO_3(\bbC)\times \SO_2(\bbC)\), realized (up to finite central subgroup) as the stabilizer of the conic in the automorphism group \(\SO_5(\bbC)\) of the three dimensional quadric. 
This is the Fano threefold number 2-29. 
\end{exa}

\subsection{The moment polytope, weight lattice and valuation cone}

Let \(X\) be a projective spherical \(G\)-variety. 
To a fixed \(G\)-linearized line bundle \(L\) on \(X\) is associated a sequence of representations of \(G\): the spaces of plurisections \(H^0(X,L^k)\) for \(k\in \bbN_{>0}\). 
Each of these spaces admit a decomposition as a (multiplicity free) direct sum of irreducible representations: 
\[ H^0(X,L^k) = \sum_{\lambda\in \Lambda_k} V_{\lambda} \]
where \(\Lambda_k\) is a (finite) subset of the set of dominant weights of \(G\).

\begin{defn}
The \emph{moment polytope} of \((X,L)\) is the set 
\(\Delta = \Delta(X,L) = \overline{\bigcup_{k\in \bbN_{>0}}\frac{1}{k}\Lambda_k }\).
\end{defn}

The name comes from the coincidence with the symplectic moment polytope of Kirwan, and it is indeed a convex polytope \cite{Brion_1999}. 
To state the K-stability criterion for spherical varieties, we need two additional pieces of information encoded by combinatorial objects. 

\begin{defn}
The \emph{weight lattice} \(M\subset X^*(B)\) is the set of weights of \(B\)-semi-invariant rational functions on \(X\). 
\end{defn}

The rank of \(M\) is called the \emph{rank} of the spherical \(G\)-variety \(X\), and coincides with the dimension of the moment polytope of any \(G\)-linearized ample line bundle \(L\) on \(X\). 
The weight lattice allows to recover \(H^0(X,L)\) from the moment polytope \(\Delta(X,L)\) and a rational \(B\)-semi-invariant section \(s\) of \(L\). 

\begin{thm}[\cite{Brion_1989}]
If \(\chi\) is the weight of \(s\), then 
\[ H^0(X,L) = \bigoplus_{\lambda\in (\chi + M)\cap \Delta} V_{\lambda} \]
\end{thm}

Let \(N=\Hom(M,\bbZ)\) be the dual lattice. 
Any \(\bbQ\)-valued valuation \(\nu\) on the field \(\bbC(X)\) of rational functions on \(X\) determines an element of \(N\otimes \bbQ\) denoted by \(\nu|_M\) by abuse of notations, defined by the values the valuation takes on \(B\)-semi-invariant rational functions. 

\begin{defn}
The \emph{valuation cone} of \(X\) is the set \(\mathcal{V}\subset N\otimes \bbQ\) of all \(\nu|_M\) where \(\nu\) runs over the set of \(G\)-invariant \(\bbQ\)-valued valuations on \(\bbC(X)\). 
\end{defn}

We refer to \cite{Delcroix-Montagard} for the explicit data of valuation cone, weight lattice for locally factorial spherical Fano varieties of dimension less or equal to four and rank less or equal to two as examples, together with the (dual of the) moment polytope for the anticanonical line bundle. 
Examples of moment polytopes for arbitrary ample line bundles may be found in \cite{Delcroix_ZAG,Delcroix_2023_RK1,Delcroix_2023_unif}. 

\section{K-stability for spherical varieties}

We first review the results of \cite{Delcroix_2023_unif} which provide a convex geometric translation of K-stability for polarized spherical varieties. 
These results generalize the case of toric varieties treated by Donaldson in \cite{Donaldson_2002}, and partial results in the direction of spherical varieties were obtained by Alexeev and Katzarkov previously \cite{Alexeev-Katzarkov_2005}.
We then review various cases where the convex geometric criterion may further be reduced to an effectively checkable criterion, providing examples of progress towards Problem~\ref{effective_YTD}. 

\subsection{Test configurations as convex PL functions}

Fix \(\chi\) in the affine span of \(\Delta\) in \(X^*(B)\otimes \bbR\), which is parallel to \(M\otimes \bbR\). 
In particular, \(-\chi+\Delta\subset M\otimes \bbR\), and an element \(v\in \mathcal{V}\subset N\otimes\bbR=(M\otimes \bbR)^*\) defines a linear form on \(M\otimes \bbR\). 
We introduce a cone \(\mathcal{C}\) of convex functions on \(-\chi+\Delta\) by setting 
\[ \mathcal{C} := \{f=\sup(c_1-v_1,\ldots,c_k-v_k):-\chi+\Delta \to \bbR \mid k\in \bbZ_{>0}, c_j\in \bbR, v_j\in \mathcal{V} \} \]
In other words, a function on \(-\chi+\Delta\) is in \(\mathcal{C}\) if it is a convex, piecewise rational affine function all of whose slopes lie in \(-\mathcal{V}\). 

Given \(f=\sup(c_1-v_1,\ldots,c_k-v_k) \in \mathcal{C}\), we let \(\nld(f)\) be the number of linearity domains of \(f\) on \(-\chi +\Delta\), that is, the number of connected components of \(M\otimes \bbR\setminus \bigcup_{1\leq i,j \leq k} \{c_i-v_i=c_j-v_j\} \) whose relative interior intersects \(\Delta\). 

\begin{thm}[\cite{Delcroix_2023_unif}]
\label{TC_classification}
Let \((X,L)\) be a spherical polarized \(G\)-variety. 
Then \(G\)-equivariant test configurations for \((X,L)\) are encoded by functions in \(\mathcal{C}\). 
Furthermore, the number of irreducible components of the central fiber \(X_0\) of a test configuration \(\mathcal{X}\) associated with the function \(f\) is equal to \(\nld(f)\). 
The test configuration is product (up to base change) if and only if the function \(f\) is affine on \(-\chi+\Delta\), and defined by an element of the linear part \(\mathcal{V}\cap -\mathcal{V}\) of the valuation cone. 
\end{thm}

In the above, the translation by \(\chi\) plays only a superficial role, but it is useful to keep in mind that, in general, \(\Delta\) does not live in \(M\otimes \bbR\) but in a translate. 

While we will not carry out the proof here, let us sketch the correspondence. 
Let \(\mathcal{X}\) be a \(G\)-equivariant test configuration for \((X,L)\), and compactify it to a family \(\bar{\mathcal{X}}\) over \(\bbP^1\) by gluing the trivial family at infinity as usual in K-stability. 
Then the normal variety \(\bar{\mathcal{X}}\) is spherical under the action of \(G\times \bbC^*\), of rank equal to the rank of \(X\) plus one. 
Up to twisting by a line bundle pulled back from the base \(\bbP^1\), \(\mathcal{L}\) extends to a \(G\times \bbC^*\)-linearized ample line bundle \(\bar{\mathcal{L}}\) on \(\bar{\mathcal{X}}\). 
Consider the associated moment polytope \(\Delta(\bar{\mathcal{X}},\bar{\mathcal{L}})\). 
Applying the results in \cite{Brion_1989}, one can show that it is cut out by a copy of \(\Delta\), walls defined by the same normals as the normals of \(\Delta\), and additional facets which are exactly described by the affine parts of the graph of a function \(f\in \mathcal{C}\). 
The converse construction follows roughly from the classification of spherical varieties by colored fans \cite{Knop_1991} and \cite{Brion_1989}. 

\subsection{The non-Archimedean Mabuchi functional and the K-stability criterion}

Let \(\mathop{d\mu}\) be the Lebesgue measure on the affine span of \(\Delta\), normalized by the translated lattice \(\chi+M\) (that is, if \((e_1,\ldots,e_r)\) generates the lattice, then the parallelepiped whose sides are given by the \(e_j\), translated by \(\chi\), has mass \(1\) for the measure \(d\mu\)). 
For any facet \(F\) of \(\Delta\), let \(\mathop{d\sigma|_F}\) be the Lebesgue measure on the affine span of \(F\), normalized by the lattice obtained by intersection with (a translate of) \(M\). 

Recall that \(\Phi^+\) denotes the set of positive roots of \(G\) with respect to the torus \(T\) and the Borel subgroup \(B\), and let \(\varpi\) be the half-sum of all positive roots. 
Let \(\kappa\) denote the Killing form of \(G\). 
Let \(\Phi_X^+\) denote the set of positive roots \(\alpha\in \Phi^+\) such that there exists \(p\in \Delta\) with \(\kappa(\alpha,p)\neq 0\).  
We consider two polynomials defined on \(X^*(B)\otimes \bbR\) by 
\[ P(p) = \prod_{\alpha\in \Phi_X^+} \frac{\kappa(\alpha,p)}{\kappa(\alpha,\varpi)} \]
and 
\[ Q(p) = \sum_{\alpha\in \Phi_X^+}\frac{\kappa(\alpha,\varpi)}{\kappa(\alpha,p)}P(p) \]
We consider the functional \(\mathcal{L}\) defined on the cone of functions \(\mathcal{C}\) by 
\[ \mathcal{L}(f) = \int_{\partial \Delta} f(p-\chi)P(p)\mathop{d\sigma}(p) - \int_{\Delta}f(p-\chi) (aP(p)-Q(p))\mathop{d\mu}(p) \]
where \(a\in \bbR\) is uniquely determined by requiring that \(\mathcal{L}\) vanishes on constant functions. 

\begin{thm}[\cite{Delcroix_2023_unif}]
For the \(G\)-equivariant test configuration \(\mathcal{X}\) for \((X,L)\) associated with the function \(f\in \mathcal{C}\), the non-Archimedean Mabuchi functional reads 
\[ \MNA(\mathcal{X}) = \mathcal{L}(f)  \]
\end{thm}

Thanks to Chi Li's result \cite{Li_2022} as applied to spherical varieties by Yuji Odaka in the appendix to \cite{Delcroix_2023_unif},
or in special cases \cite{Delcroix_2020_horo,Li-Zhou-Zhu_2018} direct approaches using \cite{Chen-Cheng_2021}, 
we obtain a criterion for the existence of cscK metrics on polarized spherical varieties. 

\begin{thm}[\cite{Delcroix_2023_unif}]
A smooth polarized spherical \(G\)-variety \((X,L)\) admits a cscK metric in \(c_1(L)\) if and only if, for any \(f\in \mathcal{C}\), one has \(\mathcal{L}(f)\geq 0\), with equality if and only if \(f\) is affine defined by an element of the linear part \(\mathcal{V}\cap -\mathcal{V}\) of the valuation cone. 
\end{thm} 

This characterization simplifies greatly the K-polystability criterion: the space of test configurations becomes a convex cone in the space of convex functions on the moment polytope, and the non-Archimedean Mabuchi functional is linear with respect to this structure. 
However, the cone of functions to consider is still infinite dimensional and one cannot check the condition directly. 
Provided a positive answer to the effective YTD conjecture (Problem~\ref{effective_YTD}), the cone of functions to consider would become a finite dimensional one, and one could check the condition effectively (at the very least, one could use numerical approximations to test K-stability). 
This is why we refer to Problem~\ref{effective_YTD} as an effective version of the YTD conjecture.

\subsection{Effective K-stability of Fano spherical varieties}

For a Fano spherical variety, with \(L=K_X^{-1}\), the criterion becomes much simpler (and generalizes previous results \cite{Wang-Zhu_2004,Podesta-Spiro_2010,Delcroix_2017}). 

\begin{thm}[\cite{Delcroix_2020}]
A smooth spherical Fano \(G\)-variety \(X\) admits a Kähler-Einstein metric if and only if the barycenter of its moment polytope \(\Delta(X,K_X^{-1})\) with respect to the measure \(P\mathop{d\mu}\) is in \(2\varpi_X + \relint(-\mathcal{V})^{\vee}\), where \(2\varpi_X=\sum_{\alpha\in \Phi_X^+}\alpha\). 
\end{thm}

We refer for example to \cite{Delcroix_2020,Delcroix_2022,Delcroix-Montagard,Lee-Park-Yoo_2021,Lee-Park-Yoo_2022,Hong-Hwang-Park_2024} for examples of explicit Fano manifolds on which this criterion is applied. 
Other effective results on K-stability of Fano manifolds include the case of complexity one \(T\)-varieties \cite{Ilten-Suss_2017}, and more recently the case of complexity one \(G\)-varieties in general \cite{Li-Li}, unifying the two approaches. 
A huge effort has also been put in the recent years in determination of K-stability for all Fano threefolds, starting with the already stunning book \cite{the_book} on general elements in each deformation family.

\subsection{Effective K-stability of rank one spherical varieties}

Quite surprisingly, we were able to show that for spherical varieties of rank one, the same phenomenon as in the Fano case occurs: it is enough to test K-stability on special test configurations. 
Furthermore, it follows then from the classification of such test configurations (Theorem~\ref{TC_classification}) that there is, up to rational multiple, only one non-trivial special equivariant test configuration in this setting. 
Finally, this non-trivial special equivariant test configuration is product if and only if there is a \(\bbC^*\)-action on \(X\) commuting with the action of \(G\), that is, if and only if the \(G\)-variety \(X\) is horospherical. 

\begin{thm}[\cite{Delcroix_2023_RK1}]
Let \((X,L)\) be a smooth, rank one polarized spherical \(G\)-variety, and let \(\ell\in \mathcal{V}\setminus \{0\}\). 
Then there is a cscK metric in \(c_1(L)\) if and only if 
\(\mathcal{L}(\ell)\geq 0\), with equality if and only if \(X\) is horospherical. 
\end{thm}

This gives the best possible answer in the cscK setting to the question raised by Chen and Guan in 2000 \cite{Guan-Chen_2000}: is it possible to find a special metric on a compact almost homogeneous manifold of cohomogeneity one?
Indeed, a polarized Kähler manifold of cohomogeneity one is a spherical rank one variety by \cite{Cupit-Foutou_2003}. 
Examples of applications of this result are provided in \cite{Delcroix_2023_RK1, Delcroix_ZAG}, in particular, we provide examples where all Kähler classes admits cscK metrics, as well as examples with both cscK and non-cscK Kähler classes. 

Note that rank one horospherical varieties are, up to blowing up the closed orbits, \(\bbP^1\)-bundles over generalized flag manifolds. 
For more general \(\bbP^1\)-bundles satisfying some geometrical assumptions, called semisimple principal \(\bbP^1\)-bundles, Apostolov, Calderbank, Gauduchon and Tønnesen-Friedmann prove in \cite{ACGTIII_2008} that they admit cscK metrics (and more generally extremal Kähler metrics) if and only if the associated so-called extremal polynomial is positive. 
The latter condition may be reinterpreted as some sort of effective YTD conjecture: it is equivalent to K-polystability for (some) real test configurations whose central fiber has at most two irreducible components. 

\subsection{Effective K-stability of toric surfaces}

To finish with results on the effective YTD conjecture, let us go back in time a bit, and review an amazing statement due to Donaldson, with some additional input of Wang and Zhou, applying to toric surfaces. 
As we said, toric varieties are spherical varieties. 
Furthermore, the rank of a toric variety is the same as its dimension, so we are considering a special subclass of rank two spherical varieties. 

\begin{thm}[\cite{Donaldson_2002,Donaldson_2009,Wang-Zhou_2014}]
A smooth polarized toric surface \((X,L)\) admits a cscK metric in \(c_1(L)\) if and only if for any two affine functions \(\ell_1\) and \(\ell_2\) on \(M\otimes \bbR\), \(\mathcal{L}(\sup(\ell_1,\ell_2))\geq 0\), with equality if and only if \((\ell_1-\ell_2)^{-1}(0)\) contains or does not meet the relative interior of \(\Delta\) (i.e. \(\sup(\ell_1,\ell_2)\) is affine on \(\Delta\)).  
\end{thm}

We should mention that the result of Wang and Zhou is also stated for extremal Kähler metrics in \cite{Wang-Zhou_2014}. 
The criterion is applied to many examples of toric surfaces in \cite{Donaldson_2002,Wang-Zhou_2011,Wang-Zhou_2014,Wei-Chen_2020}.

\bibliographystyle{alpha}
\bibliography{SKSTS}
\end{document}